\documentclass[a4paper,11pt]{article}
\usepackage[colorlinks,bookmarksopen,bookmarksnumbered,colorlinks = false]{hyperref}
\usepackage{caption}
\usepackage[a4paper, top=1.5in, total={5.5in, 8.5in}]{geometry}
\usepackage{graphicx}
\usepackage{array}
\usepackage[ruled,vlined]{algorithm2e}
\usepackage{pgfplots}
\pgfplotsset{compat=1.8}
\usepackage{amssymb}
\usepackage{amsmath}
\usepackage{fancyhdr}

\fancypagestyle{alim}{\fancyhf{}\fancyfoot[L]{\fontsize{9}{11} \selectfont  $^*$email: \texttt{T.J.Deveney@bath.ac.uk}\\ $^\dagger$email: \texttt{E.Mueller@bath.ac.uk}\\$^\ddagger$email: \texttt{T.Shardlow@bath.ac.uk}}}
\usepackage{bm}
\usepackage{float}
\newcommand\numberthis{\addtocounter{equation}{1}\tag{\theequation}}

\usepackage[latin1]{inputenc}
\usepackage{tikz}
\usepackage{multirow}
\usetikzlibrary{trees}

\usepackage{amssymb}
\usepackage[colorlinks,bookmarksopen,bookmarksnumbered]{hyperref}
\usepackage{mathtools,xparse}
\usepackage{caption}
\title{A deep surrogate approach to efficient Bayesian inversion in PDE and integral equation models}
\author{\fontsize{13}{11} \selectfont Teo Deveney$^{*}$, Eike Mueller$^{\dagger}$, Tony Shardlow$^{\ddagger}$ \\ \fontsize{10}{11} \selectfont Department of Mathematical Sciences, University of Bath, Bath, UK, BA2 7AY}

\renewcommand{\theequation}{\arabic{equation}}

\date{}
\begin{document}

\numberwithin{equation}{section}

\maketitle
\thispagestyle{alim}
\begin{abstract}
\noindent
We investigate a deep learning approach to efficiently perform Bayesian inference in partial differential equation (PDE) and integral equation models over potentially high-dimensional parameter spaces. The contributions of this paper are two-fold; the first is the introduction of a neural network approach to approximating the solutions of Fredholm and Volterra integral equations of the first and second kind. The second is the development of a new, efficient deep learning-based method for Bayesian inversion applied to problems that can be described by PDEs or integral equations. To achieve this we introduce a surrogate model, and demonstrate how this allows efficient sampling from a Bayesian posterior distribution in which the likelihood depends on the solutions of PDEs or integral equations. Our method relies on the direct approximation of parametric solutions by neural networks, without need of traditional numerical solves. This deep learning approach allows the accurate and efficient approximation of parametric solutions in significantly higher dimensions than is possible using classical discretisation schemes. Since the approximated solutions can be cheaply evaluated, the solutions of Bayesian inverse problems over large parameter spaces are efficient using Markov chain Monte Carlo. We demonstrate the performance of our method using two real-world examples; these include Bayesian inference in the PDE and integral equation case for an example from electrochemistry, and Bayesian inference of a function-valued heat-transfer parameter with applications in aviation. 
\end{abstract}

\newpage
\section{Introduction}
Deep learning methodologies have seen significant development in the past few decades. The advent of efficient optimisation algorithms has led to impressive results in many high-dimensional tasks, particularly in settings where there is an abundance of data. Recently there has been much focus on the use of  deep learning to address mathematical challenges involving differential equations. In the short time of just a few years, this has led to much progress, both in algorithmic ingenuity and theoretical understanding of neural network approaches to specific mathematical problems. For example, we now have proof that deep neural networks can overcome the curse of dimensionality in the approximation of the solutions of certain classes of PDEs \cite{jen1,jen4,jen5}, and efficient algorithms capable of generating these approximations have been identified \cite{han17,rai18,e17,sir17}. These methods have already found applications in a range of areas such as fluid dynamics \cite{rai188} and financial mathematics \cite{zh19,be18}. Further exploration of these methods is well underway, with applications including optimal control problems using the Hamilton-Jacobi-Bellman equation \cite{hjb1,hjb2}, and PDE-based regularisation in regression tasks for which the governing laws are partially known \cite{physreg}. In forward uncertainty quantification, surrogate models and probabilistic generative models have been implemented to capture response uncertainty of PDEs with random terms \cite{sur11},  \cite{sur111}. 
\\ \\
While deep learning approximations for PDE-based problems have received a lot of attention, integral equations \cite{ra_int} are less well investigated, despite also having widespread applications in areas such as radiative transfer, viscoelasticity, and electrochemistry. This paper extends previous work on the solution of PDE problems to integral equations, by introducing a new method for solving integral equations using neural networks. Furthermore, existing approaches for PDEs focus mostly on the solution of these equations and forward uncertainty quantification. The second novelty of this paper is the application of deep learning to Bayesian inverse problems for which the likelihood function depends on the solution of a PDE or integral equation. Our method is based on the direct approximation to the solution of the parametric forward problem by deep surrogate models. In both the PDE and integral equation case, our surrogates are analytic, easily differentiable, and allow for rapid approximate evaluations of the solution for different model parameters without need of further numerical solves. This makes the full library of existing parametric statistics techniques accessible for these physics-based problems without modification. Using these surrogates we apply Markov chain Monte Carlo (MCMC) to rapidly sample from the Bayesian posterior distribution of uncertain parameters given data and a prior distribution. This is a novel application of parametric forward solvers based on deep neural networks, which achieves a significant speed-up when compared to a more traditional finite difference MCMC approach, without sacrificing accuracy. Since deep neural networks overcome the curse of dimensionality, our approach is efficient and accurate for high-dimensional problems or problems with a large parameter space, for which the repeated approximation of the forward problem becomes computationally intractable using more traditional discretisation schemes. 
\\ \\ 
Multiple alternative approaches have been devised to improve the efficiency of MCMC methods for inverse problems, such as the use of inexpensive coarse scale models as pre-conditioners for the fine scale model \cite{ef06,ef05}, and multilevel techniques \cite{mlmcmc}. These methods can reduce the number of expensive numerical solves of the fine scale problem that result in rejected proposals in cases where the solutions at the coarser levels are closely correlated to the true solution, however a large number of fine scale solves are still required for sufficient new proposals to be accepted. Surrogates based on approximations to the parametric forward problem have also been employed. A common approach is to use a generalised polynomial chaos representation constructed by a stochastic Galerkin projection \cite{polcha1,stogal}. In this scheme a orthogonal polynomial decomposition is used to approximate the parametric solution, however this becomes intractible in higher dimensions due to the curse of dimensionality of the polynomial approximation to the solution. For this reason, the approach has been restricted to applications with few model parameters \cite{ry19,he14}. Collocation based surrogates have been applied in cases where the number of parameters is larger (though typically less than 10). In this case, more expensive traditional discretisation schemes are employed to generate an ensemble of solutions for a range of parameter values, and these are interpolated to construct the surrogate \cite{col1,polcha}. This approach requires a large number of expensive solves, and relies on interpolants which may not satisfy the required equations in complicated regions of the domain \cite{sur111}. While all of these methods have all been successful in accelerating MCMC for various applications, each is restricted either by the dimensionality of the problems that they can be applied to, or the accuracy of the surrogates that they construct. Our approach overcomes these issues by using a deep learning approximation that is accurate in high dimensions and rapid to evaluate.
\\ \\
The key contributions of this paper are as follows:

\begin{itemize}
    \item We introduce a new deep learning approach for the solution of parametric Fredholm and Volterra integral equations of the first and second kind. By capturing the dependence on model parameters, our method provides an efficient approximation of the parametric forward problem for high-dimensional integral equation parameter spaces.
 
\item We describe how neural network surrogate models for PDE and integral equation based problems can be used to rapidly sample from the posterior distribution in a Bayesian setting. In particular, we argue that our approach is efficient even if the problem is high-dimensional or depends on a large number of parameters.
 
\item By considering a specific example from electrochemistry, we demonstrate that our deep surrogate integral-equation solvers can be used to infer parameters in a Bayesian setting with an accuracy that is comparable to current PDE- based deep learning methods.
 
\item To demonstrate that the method can efficiently perform Bayesian inversion in a high-dimensional parameter space, we apply our deep surrogate solver to infer a parametric function in the PDE for heat transfer in aircraft turbines.
\end{itemize}
This paper is organised as follows: In Section 2, we describe the general form of the PDE and integral equation models, and define the deterministic and parametric forward problem that will be considered. Here we also introduce the statistical model that form the basis of the Bayesian inverse problem \cite{stuart}. In Section 3, we review the physics-informed neural network \cite{pidl1} and the deep Galerkin method \cite{sir17} for solving the PDE forward problem, and describe an extension of these methods to the accurate approximation of the parametric forward problem by a deep surrogate model. In Section 4, we introduce a novel deep learning method for the solution of integral equations and outline its extension to parametric integral equation problems using a deep surrogate model. In Section 5, we explain how to apply a MCMC scheme that uses deep surrogate model evaluations to significantly reduce the computational cost of sampling from the posterior distribution. Finally in Section 6 we demonstrate how the deep surrogate method can accelerate Bayesian inference by applying our method to two problems with real applications in electrochemistry and heat transfer in aircraft turbines. These examples cover PDE-based problems and integral equations, with inferred parameters including scalar values and functions. Python code is available to reproduce all of the presented examples \cite{zenodo_doi}.

\section{Problem specification}
We begin by defining the general form of the PDE and integral equation problems that will be investigated in this paper. This includes a description of the deterministic and parametric forward problems for these models, as well as the inverse problem corresponding to these models.
\subsection{Forward problem}
First consider the parameterised scalar-valued PDE,
\begin{align}
    \qquad \qquad \mathcal{N}(\bm{x},u(\bm{x};\bm{\theta}) ;\bm{\theta}_N)  &= h(\bm{x};\bm{\theta}_h),  \qquad \bm{x} \in \Omega,\ \bm{\theta} \in \Theta,  \label{pde}
\end{align}     
with boundary conditions defined by
\begin{align}
    \qquad \qquad \qquad u(\bm{x};\bm{\theta}) &= b(\bm{x};\bm{\theta}_b), \qquad \qquad  \bm{x} \in   \partial \Omega,\ \bm{\theta} \in \Theta. \label{bc}
\end{align}
Here $\Theta \subset \mathbb{R}^{^p}$ denotes the parameter space and $\bm{\theta} = (\bm{\theta}_N,\bm{\theta}_h,\bm{\theta}_b) \in \Theta$ is a parameter vector. The input space is denoted $\Omega \subset \mathbb{R}^d$, where $d$ is the dimension of the PDE domain $\Omega$. We assume that $\mathcal{N}$ is a known (possibly non-linear) differential operator parameterised by $\bm{\theta}_N$, and $h,b$ are known functions parameterised by $\bm{\theta}_h$ and $\bm{\theta}_b$. If $\bm{\theta}$ is fixed this is a deterministic PDE, and the deterministic forward problem is to approximate the $d$-dimensional function $u_\theta(\bm{x}):\Omega \to \mathbb{R}$ satisfying (\ref{pde}, \ref{bc}) given $\bm{\theta}$. Here we have used $u_\theta(\bm{x})$ in place of $u(\bm{x};\bm{\theta})$ as the solution to (\ref{pde},\ref{bc}) to signify that the parameters $\bm{\theta}$ are fixed, and therefore not inputs of the solution. This convention is extended to all functions and operators throughout this work, whereby $\bm{\theta}$ is placed into their subscript when it is fixed. The parametric forward problem is to produce a function capable of returning the solution corresponding to any parameters $\bm{\theta}\in \Theta$, that is to approximate the $d+p$-dimensional function $u(\bm{x};\bm{\theta}):\Omega \times \Theta \to \mathbb{R}$ satisfying (\ref{pde}, \ref{bc}). 
\\ \\
In this work, we also consider non-homogenous integral equations of the first kind
\begin{align}
    \qquad 0 = v(x;\bm{\theta}_v) + \int_a^{b(x)} k(x,y;\bm{\theta}_k)u(y;\bm{\theta}) dy, \qquad x \in \Omega, \ \bm{\theta} \in \Theta, \label{firstkind}
\end{align}
 and the second kind
 \begin{align}
    \qquad u(x;\bm{\theta}) = v(x;\bm{\theta}_v) + \int_a^{b(x)} k(x,y;\bm{\theta}_k)u(y;\bm{\theta}) dy,\qquad  x \in \Omega, \ \bm{\theta} \in \Theta. \label{inteq}
\end{align}
Again $\Theta \subset \mathbb{R}^{^p}$ denotes the parameter space and $\bm{\theta} = (\bm{\theta}_v,\bm{\theta}_k) \in \Theta$ is a parameter vector. The input space is $\Omega = [a,b^*]$ for some $a<b^*$. We assume $v(x;\bm{\theta}_v)$ is a given function parameterised by $\bm{\theta}_v$, and $k(x,y;\bm{\theta}_k)$ is a known kernel parameterised by $\bm{\theta}_k$. The deterministic forward problem is to approximate the function $u_\theta(x):\Omega \to \mathbb{R}$ satisfying the integral equation for a fixed $\bm{\theta}$. The parametric forward problem is to approximate $p+1$-dimensional function $u(x;\bm{\theta}):\Omega \times \Theta \to \mathbb{R}$ satisfying the equation. Our deep learning algorithm can solve the parametric forward problems (\ref{firstkind}) and (\ref{inteq}) for both Fredholm equations ($b(x) = b^*$), and Volterra equations ($b(x) = x$). 

\subsection{Inverse problem}
The goal of the inverse problem for both the PDE and integral equation is to estimate $\bm{\theta}$ given some data samples $(\hat{\bm{x}}_i,\hat{z}_i)_{i=1:M}$. This data consists of input data $\hat{\bm{x}}_i \in \Omega$ and measured responses $\hat{z_i} \in \mathbb{R}$ corresponding to those inputs. We model the data as the solution of the equation with an additional error term
\begin{align}
    \qquad \qquad \qquad \qquad \hat{z}_i = u(\hat{\bm{x}}_i;\bm{\theta}) + \epsilon_i, \qquad \qquad i=1,\dots,M. \label{model}
\end{align}
This is similar to the common scenario in parametric statistics
\begin{align}
    \qquad \qquad \qquad \qquad \hat{z}_i = f(\hat{\bm{x}}_i;\bm{\theta}) + \epsilon_i, \qquad \qquad i=1,\dots,M. \label{stat}
\end{align}
In statistical applications $f(\bm{x}_i;\bm{\theta})$ is typically a carefully constructed parametric function, designed using specialist knowledge to take input variables $\bm{x}_i$ and produce outputs $z_i$ which model the mean behaviour of the system. In our case this function is the parametric solution $u(\hat{\bm{x}}_i;\bm{\theta})$, therefore it behaves according to the laws imposed by the PDE or integral equation. The last term is in \eqref{model} a random variable which accounts for variation of the data around the mean. A standard choice for the $\epsilon_i$, which we will adopt, is to assume unbiased i.i.d. Gaussian deviations of the form $\epsilon_i = N(0,\sigma^2)$, where $\sigma^2$ is an unknown parameter. 
\\ \\
The objective in the inverse problem is the same as a statistical regression problem, that is to infer the unknown parameters $\{\bm{\theta},\sigma^2\}$ such that the model describes the data. We note that using our deep surrogate methodology, the inverse problem can be solved with any parametric regression scheme from statistics without modification, since the deep surrogate model provides a functional representation of the solution $u(x,\bm{\theta})$ that can be rapidly evaluated and differentiated with respect to the parameters $\bm{\theta}$. In this work we focus on applying MCMC methods from Bayesian parametric regression to solve Bayesian inverse problems for which the forward problem can be described by a PDE and integral equation. This provides a natural form of inverse uncertainty quantification in the form of a posterior probability distribution over all parameters of interest, but has proven too expensive to be implemented fully in the majority of cases with standard numerical discretisation techniques such as finite differences or quadrature methods.

\section{Review of deep learning techniques for PDEs}
Neural networks are capable of accurately approximating high-dimensional functions using significantly fewer parameters than traditional function approximations based on linear basis expansions. For multiple classes of PDEs it can be proven that the number of parameters required to approximate the solution to a fixed accuracy grows at most polynomially with the dimension of the problem \cite{jen1,jen4,jen5}, in contrast to exponentially for techniques like finite differences and finite elements. Because of this, deep learning techniques have shown great success in the solutions of PDE forward problems with dimensions of up to $d=200$ \cite{han17,rai18,e17,sir17}, which would have been intractable to solve with traditional discretisation approaches. Here we briefly review the physics informed neural network and deep Galerkin methods of approximating deterministic PDE solutions by neural networks. We also describe an extension of these methods to approximate the solution to the parametric forward problem which will be used in our Bayesian approach to inverse problems.
\subsection{Deterministic forward problem}
We first consider the neural network approximation of the deterministic PDE forward problem defined by $(\ref{pde},  \ref{bc})$. Recent work in deep learning for PDEs \cite{pidl1,sir17} has uncovered a simple method designed to solve such equations over a closed domain $\Omega$. Let $\pi^d$, and $\pi^b$ be measures with supports supp$(\pi^d)=\Omega$ and supp$(\pi^b)=\partial \Omega$. Under the assumption that the PDE is well posed, the approach in \cite{pidl1,sir17} exploits the fact that for a fixed $\bm{\theta}\in \Theta$ and $\nu_1,\nu_2>0$ the quantity
\begin{align*}
        \nu_1\|\mathcal{N}_\theta(u_\theta(\bm{x}),\bm{x}) - h_\theta(\bm{x})&\|^2_{L_2(\Omega,\pi^d)} + \nu_2\|u_\theta(\bm{x}) - b_\theta(\bm{x})\|^2_{L_2(\partial \Omega,\pi^b)}. \numberthis \label{lossl2}
\end{align*}
is non-negative, and zero only when $u_\theta(\bm{x})$ is the PDE solution. Here the subscript $L_2(D,\mu)$ denotes the $L_2$-norm over domain $D$ with respect to the measure $\mu$ 
\begin{align}
    \|f(\bm{x})\|^2_{L_2(D,\mu)} = \int_D|f(\bm{x})|^2d\mu.
\end{align}
The universal approximation theorem \cite{cyb89} states that neural networks are dense in the space of continuous functions. It follows that if a solution exists and is continuous, then it can be approximated arbitrarily well by a neural network. The deep learning methods described in \cite{pidl1,sir17} use neural networks to approximate the solutions of PDEs. Denoting the neural network approximation as $\hat{u}_\theta(\bm{x}) \approx u_{\theta}(\bm{x})$, these algorithms take inspiration from \eqref{lossl2} by using some variant of the gradient descent scheme to update the neural network parameters to minimise the loss function
\begin{align*}
    \frac{\nu_1}{N}\sum_{n=1}^N (\mathcal{N}_\theta(\hat{u}_\theta(\bm{x}^n),\bm{x}^n) - h_\theta(\bm{x}^n)&)^2 + \frac{\nu_2}{J}\sum_{j=1}^J(\hat{u}_\theta(\bm{y}^j) - b_\theta(\bm{y}^j))^2. \numberthis \label{loss}
\end{align*}
Here $\bm{x}^n \in \Omega,\ \bm{y}^j \in  \partial \Omega $ are collocation points. If these points are distributed according to the measures $\pi^d, \pi^b$ then \eqref{loss} is a Monte Carlo estimate of \eqref{lossl2}, and the network is trained to approximate the solution to $(\ref{pde},\ref{bc})$ by minimising this estimate with respect to the network parameters. The coefficients $\nu_1,\nu_2$ are interpreted as positive weights which can be adjusted to reflect the relative importance of the terms in \eqref{loss}. The minimisation algorithm depends on gradients, which are computed via back-propagation using automatic differentiation. Automatic differentiation is applied to compute both the descent direction of the network parameters, and to apply the differential operator $\mathcal{N}$ to the neural network at the collocation points. These computations can be routinely implemented using automatic differentiation software such as TensorFlow \cite{tensorflow}. 
\\ \\
If the collocation points are constant throughout training, we arrive at the physics-informed neural network method in \cite{pidl1}. This method is based on a fixed mesh, which requires a sufficiently fine resolution in order to yield accurate results, since the neural network must interpolate the solution between training points. In higher-dimensional spaces, this becomes prohibitively expensive since the number of collocation points required to saturate $\Omega$ grows exponentially with the PDE dimension $d$. 
\\ \\
In contrast, if the collocation points are randomly sampled from $\pi^d,\pi^b$ before each gradient descent iteration, the minimisation of \eqref{lossl2} results in the deep Galerkin method (DGM) described in \cite{sir17}. This mesh-free approach has proven efficient at accurately approximating high-dimensional solutions, since it allows for substantially fewer collocation points to be used at each iteration. By re-sampling before each iteration, it relies on a mini-batch stochastic gradient descent approach in order to adequately cover the domain. In \cite{sir17} this approach was used to solve a PDE with dimension $d=200$. An outline of this algorithm is: \\ \\
\begin{algorithm}[H]
\SetAlgoLined
    1.  Initialise a neural network approximation of the solution $\hat{u}_\theta(\bm{x}):\mathbb{R}^{d} \to \mathbb{R}$ \\
    2.  Draw i.i.d. random samples $(\bm{x}^n)_{n=1:N} \sim \pi^d, (\bm{y}^j)_{j=1:J} \sim \pi^b$ \\
    3.  Update the neural network parameters by taking \\ \ \ \ \ one gradient descent step to reduce the loss function \eqref{loss} \\
    4.  Repeat 2-3 until convergence
\caption{DGM solver for the deterministic forward PDE problem}
\end{algorithm}
\ \\
Upon convergence, this approach will approximate the solution to the PDE for a fixed choice of parameters $\bm{\theta}$. In our implementations we use the relative change in the moving average of the estimated loss function \eqref{loss} as convergence criteria.

\subsection{Parametric forward problem}
In order to adapt the deep Galerkin approach to the parametric forward problem we must introduce an additional measure $\pi^\theta$ with support supp($\pi^\theta$) = $\Theta$, and use this to extend the loss function \eqref{lossl2} so that the minimisation also spans the parameter space. The extended loss function is
\begin{align*}
        \nu_1\|\mathcal{N}(u(\bm{x};\bm{\theta}),\bm{x};\bm{\theta}_N) - h(\bm{x};\bm{\theta}_h)\|^2_{L_2(\Omega\times\Theta,\pi^d \otimes\pi^\theta)}&\\
        + \nu_2\|u(\bm{x};\bm{\theta})& - b(\bm{x};\bm{\theta}_b)\|^2_{L_2(\partial \Omega\times\Theta,\pi^b\otimes\pi^\theta)}. \numberthis \label{lossl3}
\end{align*}
In this extended loss function the $L_2$-norm is taken over the joint domain $\Omega\times\Theta$ with respect to the product measure $(\pi^d \otimes\pi^\theta)(\bm{x},\bm{\theta}) = \pi^d(\bm{x}) \pi^\theta(\bm{\theta})$. Denoting the neural network approximation of the parametric solution by $\hat{u}(\bm{x};\bm{\theta})$, we now apply stochastic gradient descent to minimise \eqref{lossl3}, by sampling $\bm{x}^n$ from $\pi^d$, $\bm{y}^j$ from $\pi^b$, and $\bm{\theta}^n,\bm{\phi}^j$ from $\pi^\theta$ at each iteration and minimising the Monte Carlo estimate given by
\begin{align*}
    \frac{\nu_1}{N}\sum_{n=1}^N (\mathcal{N}(\hat{u}(\bm{x}^n;\bm{\theta}^n),\bm{x}^n;\bm{\theta}^n_A) - h(\bm{x}^n;\bm{\theta}^n_h)&)^2 + \frac{\nu_2}{J}\sum_{j=1}^J(\hat{u}(\bm{y}^j;\bm{\phi}^j) - b(\bm{y}^j;\bm{\phi}^j_b))^2. \numberthis \label{loss8}
\end{align*} 
The algorithm for the parametric forward problem is therefore: \\ \\
\begin{algorithm}[H]
\SetAlgoLined
    1.  Initialise a neural network approximation of the parametric solution, \\ \ \ \  $\hat{u}(\bm{x};\bm{\theta}):\mathbb{R}^{d+p} \to \mathbb{R}$ \\
    2.  Draw i.i.d. samples $(\bm{x}^n)_{n=1:N} \sim \pi^d, (\bm{y}^j)_{j=1:J} \sim \pi^b$, \\ \ \ \  $(\bm{\theta_N}^n,\bm{\theta_h}^n,\bm{\theta_b}^n)_{n=1:N} \sim \pi^\theta,  (\bm{\phi_N}^j,\bm{\phi_h}^j,\bm{\phi_b}^j)_{j=1:J} \sim \pi^\theta$ \\
    3.  Update the neural network parameters by taking one gradient descent \\ \ \ \ step to reduce loss function \eqref{loss8} \\
    4.  Repeat 2-3 until convergence \\
\caption{DL solver for the parametric forward PDE problem}
\end{algorithm}
\ \\
We call the resulting neural network $\hat{u}(\bm{x}, \bm{\theta}) \approx u(\bm{x}, \bm{\theta})$ a \textit{deep surrogate model} for the parametric PDE. This method of approximating parametric solutions is efficient even over high dimensional parameter spaces where standard discretisations would exhibit exponential growth in computational requirement. In Section 6 we implement this approach to a heat transfer PDE applied in aviation, and an electrochemistry PDE problem. The results demonstrate the accuracy and efficiency of the deep surrogate model for rapid solution evaluations over a range of parameter values.
%\\ \\
%with an equivalent integral equation formulation as well as the deep learning integral equation solver that we introduce in Section 4.  

\section{A deep learning approach for integral equations}
We now extend the deep learning approaches described in the previous section to the solution of integral equations, which so far have not been discussed in the literature. Our novel methodology allows the computation of approximate solutions to the parametric forward problems for these models. 
\subsection{Deterministic forward problem}
Motivated by (\ref{firstkind}) and (\ref{inteq}), we begin by defining the function
\begin{align}
    w_\theta(x,y) := \int_a^{y}k_\theta(x,\gamma)u_\theta(\gamma)d\gamma. \label{ww}
\end{align}
Using this, for fixed parameters $\bm{\theta}$, the non-homogeneous integral equation of the second kind \eqref{inteq} can be written as
\begin{align}
    u_\theta(x) = v_\theta(x) + w_\theta(x,b(x)), \qquad x \in [a,b^*]. \label{reinteq}
\end{align}
We seek to solve \eqref{reinteq} by approximating the solution $u_\theta(x)$ by a neural network and minimising the residual norm of \eqref{reinteq} in a similar fashion to the PDE based problem discussed in Section 3. For integral equations, an added complication arises since $w(x,b(x);\bm{\theta}_k)$ is an integral involving the unknown solution $u_\theta(x)$. Differentiating \eqref{ww} with respect to $y$ shows that $w_\theta(x,y)$ solves the initial value problem
\begin{align*}
    \qquad \frac{\partial{w}_\theta}{\partial{y}}(x,y) &= k_\theta(x,y)u_\theta(y), \qquad x \in [a,b^*], \ y\in[a,b(x)] \\
    w_\theta(x,a) &= 0. \numberthis \label{integrator}
\end{align*}
%Hence if $u_\theta(x), w_\theta(x,y)$ satisfy both \eqref{reinteq} and \eqref{integrator}, then $u_\theta(x)$ is a solution to the integral equation. Proceeding similarly to the PDE based problem, we define the domain $\Omega=\{x,y:x\in[a,b^*],y\in[a,b(x)]\}$, and let $\pi$ be a measure with supp$(\pi)=\Omega$, and $\pi^x$ be the marginal supp$(\pi^x)=[a,b^*]$, then $u_\theta(x)$ solves \eqref{inteq} if $u_\theta(x)$ and $w_\theta(x,y)$ are functions such that the following is zero
%\begin{align*}
%    \nu_1\bigg\|\frac{\partial w_\theta}{\partial y}(x,y) - k(x,y|&\bm{\theta}_k)u_\theta(y)\bigg\|_{L_2(\Omega,\pi)}^2  + \nu_2\|w_\theta(x,a|\bm{\theta}_k)\|_{L_2([a,b^*],\pi^x)}^2 \\ 
%    &+ \nu_3\|u_\theta(x) - v(x|\bm{\theta}_v) - w_\theta(x,b(x)|\bm{\theta}_k)\|_{L_2([a,b^*],\pi^x)}^2 \label{intmod} \numberthis
%\end{align*}
Hence if $u_\theta(x), w_\theta(x,y)$ satisfy both \eqref{reinteq} and \eqref{integrator}, then $u_\theta(x)$ is a solution to the integral equation. Proceeding similarly to the PDE-based problem, we define the domain $\Omega=\{x,y:x\in[a,b^*],y\in[a,b(x)]\}$ and let $\pi^d$ be a measure such that supp$(\pi^d)=\Omega$. Assuming the integral equation \eqref{inteq} is well posed, then for $\nu_1,\nu_2,\nu_3>0$, $u_\theta(x)$ solves \eqref{inteq} if $u_\theta(x)$ and $w_\theta(x,y)$ are functions such that the following quantity is zero
\begin{align*}
    \nu_1\bigg\|\frac{\partial w_\theta}{\partial y}(x,y) - k_\theta(x,y&)u_\theta(y)\bigg\|_{L_2(\Omega,\pi^d)}^2  + \nu_2\|w_\theta(x,a)\|_{L_2(\Omega,\pi^d)}^2 \\ 
    &+ \nu_3\|u_\theta(x) - v_\theta(x) - w_\theta(x,b(x))\|_{L_2(\Omega,\pi^d)}^2. \label{intmod} \numberthis
\end{align*}
Here the first two terms being equal to zero ensures that \eqref{integrator} is satisfied, while the third term accounts for \eqref{reinteq}. This motivates the use of two neural networks to approximate the functions $u_{\theta}(x), w_{\theta}(x,y)$. We use  $\hat{u}_\theta(x) \approx u_{\theta}(x)$ to approximate the solution, and an \textit{integrator network} $\hat{w}_\theta(x,y) \approx w_{\theta}(x,y)$ to approximate the integral term. These networks are trained to approximately satisfy  \eqref{reinteq} and \eqref{integrator} by minimising a loss function of the form
\begin{align*}
    \frac{1}{N}\sum_{n=1}^N\bigg[\nu_1\bigg(\frac{\partial \hat{w}_\theta}{\partial y}(x^n,y^n)& - k_\theta(x^n,y^n)\hat{u}_\theta(y^n)\bigg)^2  + \nu_2\hat{w}_\theta(x^n,a)^2 \\ 
    &+ \nu_3(\hat{u}_\theta(x^n) - v_\theta(x^n) - \hat{w}_\theta(x^n,b(x^n)))^2\bigg]  \label{inteqloss} \numberthis
\end{align*}
with respect to the neural network parameters. Similarly to the PDE-based problem, \eqref{inteqloss} is an unbiased Monte Carlo estimate of \eqref{intmod} if the collocation points $(x^n,y^n) \in \Omega$ are distributed according to $\pi^d$. A mini-batch stochastic gradient descent algorithm to train a neural network approximation of the solution to integral equations of the second kind is then: \\ \\ 
\begin{algorithm}[H]
\SetAlgoLined
    1.  Initialise a neural network approximation of the solution, $\hat{u}_\theta(x):\mathbb{R} \to \mathbb{R}$ \\
    2.  Initialise an integrator network, $\hat{w}_\theta(x,y):\mathbb{R}^2 \to \mathbb{R}$ \\
    3.  Draw i.i.d. samples $(x^n,y^n)_{n=1:N} \sim \pi^d$ \\
    4.  Jointly update both neural networks by taking one gradient descent\\ \ \ \ \   step to reduce the loss function \eqref{inteqloss}\\
    5.  Repeat 3-4 until convergence
\caption{DL solver for deterministic integral equation forward problem}
\end{algorithm}
\ \\
The algorithm for solving integral equations of the first kind in \eqref{firstkind} is similar, though in this case \eqref{reinteq} is replaced by
\begin{align}
    0 = v_\theta(x) + w_\theta(x,b(x)), 
\end{align}
and the corresponding loss function becomes
\begin{align*}
    \nu_1\bigg\|\frac{\partial w_\theta}{\partial y}(x,y) - k_\theta(x,y&)u_\theta(y)\bigg\|_{L_2(\Omega,\pi^d)}^2  + \nu_2\|w_\theta(x,a)\|_{L_2(\Omega,\pi^d)}^2 \\ 
    &+ \nu_3\|v_\theta(x) + w_\theta(x,b(x))\|_{L_2(\Omega,\pi^d)}^2. \label{lossint} \numberthis
\end{align*}
Again, a Monte Carlo estimate can be constructed by using suitable random collocation points, and this can be reduced using stochastic gradient descent with respect to the parameters of a neural network to approximate the solution to the forward problem. 
\subsection{Parametric forward problem}
Following the same approach as in Section 3.2, we can also extend the integral equation solver to the solution of the parametric integral equation forward problem. For this, we introduce a measure $\pi^\theta$ with support supp($\pi^\theta$) = $\Theta$, extend the loss function \eqref{intmod} over this measure. For integral equations of the second kind the resulting loss function is
\begin{align*}
    \nu_1\bigg\|\frac{\partial w}{\partial y}(x,y;\bm{\theta}) - k(x,y;&\bm{\theta}_k)u(y;\bm{\theta})\bigg\|_{L_2(\Omega\times\Theta,\pi^d \otimes\pi^\theta)}^2  + \nu_2\|w(x,a;\bm{\theta})\|_{L_2(\Omega\times\Theta,\pi^d \otimes\pi^\theta)}^2 \\ 
    &+ \nu_3\|u(x;\bm{\theta}) - v(x;\bm{\theta}_v) - w(x,b(x);\bm{\theta})\|_{L_2(\Omega\times\Theta,\pi^d \otimes\pi^\theta)}^2. \label{intmod_par} \numberthis
\end{align*}
A neural network approximation to the solution $\hat{u}(x;\bm{\theta}):\mathbb{R}^{p+1}\to \mathbb{R}$ can then be substituted into this expression and the corresponding Monte Carlo estimate minimised with respect to the neural network parameters using stochastic gradient descent. Like in the deterministic case, parametric integral equations of the first kind can also be solved by extending \eqref{lossint} analogously. A demonstration of the algorithm applied in this case is given in Section 6, where it is compared to an accurate midpoint quadrature solution and a deep neural network approximation to an equivalent PDE formulation of the problem.

\section{Deep surrogate approach for Bayesian inference}
In the previous sections, we argued that neural networks allow the efficient approximate solution of high-dimensional parameterised problems that can be formulated as PDEs or integral equations. We now show how this can be applied to Bayesian inference. Our approach uses a deep surrogate model to approximate the solution to the parametric forward problem  $u(\bm{x};\bm{\theta})$ over the joint domain $\Omega \times \Theta$. Since the deep surrogate model is efficiently evaluated for any input values $\bm{x}\in \Omega$ and model parameters $\bm{\theta} \in \Theta$, it is inserted into a Markov chain Monte Carlo (MCMC) algorithm to quickly draw samples from the posterior distribution. We focus on the Metropolis--Hastings sampler for brevity. However a significant advantage is that our parametric approximation is easily differentiable with respect to the model parameters $\bm{\theta}$ using automatic differentiation. This allows the possibility to apply more efficient gradient based sampling methods such as Hamiltonion Monte Carlo \cite{betancourt2017conceptual}. 

\subsection{Bayesian inference}
Bayesian inference is a means of inferring a distribution over a set of parameters $\bm{\tilde{\theta}}$ by conditioning on an observed dataset $\bm{\tilde{z}}$. In this approach the posterior probability of the parameters is calculated based on a prior distribution, a likelihood function, and some observed data. The Bayesian approach allows us to make inferences about model parameters through analysis of the posterior distribution
\begin{equation}\label{bayesintro}
 p(\bm{\tilde{\theta}}|\bm{\tilde{z}}) = \frac{p(\bm{\tilde{z}}|\,\bm{\tilde{\theta}})\; p(\bm{\tilde{\theta}})}{\int p(\bm{\tilde{z}}|\,\bm{\tilde{\theta}})\; p(\bm{\tilde{\theta}}) d\bm{\tilde{\theta}}}  . 
\end{equation}
Here the prior distribution, with density $p(\bm{\tilde{\theta}})$, is chosen to describe ones existing knowledge about the parameters prior to data being observed. The likelihood, $p(\bm{\tilde{z}}|\bm{\tilde{\theta}})$ is the conditional density of the data $\bm{\tilde{z}}$ given parameters $\bm{\tilde{\theta}}$ according to a statistical model. If we assume that measurement errors on the data are normally distributed with variance $\sigma^2$ as in \eqref{model}, then the unknown parameters are $\bm{\tilde{\theta}} = \{\bm{\theta},\sigma^2\}$, and the likelihood function is 
\begin{align}
    p(\bm{\tilde{z}}|\bm{\tilde{\theta}}) = \frac{1}{{\left(2\pi\sigma^2\right)^{M/2} }}\text{exp}\left( - \frac{1}{2\sigma ^2 }\sum_{i=1}^M( {\hat{z}_i - u(\hat{\bm{x}}_i;\bm{\theta}) })^2  \right). \label{17}
\end{align}
The product of the prior and likelihood is proportional to the posterior density, $p(\bm{\tilde{\theta}}|\bm{\tilde{z}})$. This is the joint distribution of the parameters $\bm{\tilde{\theta}}$ conditional on the data $\bm{\tilde{z}}$, and represents an updated belief about the about the distribution of the unknown parameters $\bm{\tilde{\theta}}$ based on the data. From this distribution estimates of the parameter values, as well as any uncertainties and dependencies can be derived.
\\ \\
In general the exact functional form of the posterior cannot be directly computed and thus must be approximated. A popular method in statistics and data science is to apply a MCMC scheme to sample from the posterior distribution over the parameters $\bm{\tilde{\theta}}$ \cite{smith1993bayesian}. One of the simplest and most widely used MCMC schemes is the Metropolis--Hastings algorithm \cite{metropolis1953equation}. In this scheme a Markov chain is constructed by proposing successive states from a proposal distribution $q(\bm{\tilde{\theta}}_{prop}|\bm{\tilde{\theta}})$ which can depend on the current state $\bm{\tilde{\theta}}$. In this scheme it can be shown that if the proposals are accepted with probability
\begin{align}
A(\bm{\tilde{\theta}}_{prop},\bm{\tilde{\theta}}) = \text{min}\left(1, \frac{q({\bm{\tilde{\theta}}}|{\bm{\tilde{\theta}}_{prop}}) p(\bm{\tilde{z}}|\,\bm{\tilde{\theta}}_{prop})p(\bm{\tilde{\theta}}_{prop})}{q({\bm{\tilde{\theta}}_{prop}}|{\bm{\tilde{\theta}}})p(\bm{\tilde{z}}|\,\bm{\tilde{\theta}})p(\bm{\tilde{\theta}})}\right), \label{accept}
\end{align}
then $p(\bm{\tilde{\theta}}|\bm{\tilde{z}})$ is the stationary distribution of the Markov chain. Assuming mild conditions on $q(\cdot|\cdot)$ to ensure ergodicity (a sufficient condition is $q(\bm{\tilde{\theta}}|\bm{\tilde{\theta}}_{prop})>0$ $ \forall \ \bm{\tilde{\theta}},\bm{\tilde{\theta}}_{prop} \in \Theta$),  the states visited by the Markov chain form a sample from the posterior distribution. The Metropolis--Hastings algorithm can be written: \\ \\
\begin{algorithm}[H]
    \SetAlgoLined
    Choose initial $\bm{\tilde{\theta}}_0 \in \Theta$\; 
    \For{$i = 0,1,\dots,N_{samples}$}{ 
         Propose $\bm{\tilde{\theta}}_{prop} \sim q({\bm{\tilde{\theta}}_{prop}}|{\bm{\tilde{\theta}}_i})$\;
         Sample $u \sim U(0,1)$\;
         \eIf{$u \leq A(\bm{\tilde{\theta}}_{prop},\bm{\tilde{\theta}}_i)$}{
            set $\bm{\tilde{\theta}}_{i+1} := \bm{\tilde{\theta}}_{prop}$ (accept)\;
          }{
            set $\bm{\tilde{\theta}}_{i+1} := \bm{\tilde{\theta}}_i$ (reject)\;
          } 
    }
    Return $\bm{\tilde{\theta}} = (\bm{\tilde{\theta}}_0,\bm{\tilde{\theta}}_1,\dots,\bm{\tilde{\theta}}_{N_{samples}})$
    \caption{Metropolis--Hastings}
\end{algorithm} 
\ \\
The empirical distribution of samples generated in this way converges in distribution to the true posterior distribution \eqref{bayesintro} as the number of samples increases. To ensure a representative sample in practice, some proportion of the samples are removed from the start of the sample due to the 'burn-in' period. These represent the time taken for the Markov chain to transition towards the typical set of the posterior distribution and are typically a poor representation of the true posterior. Since subsequent samples in the Markov chain are correlated, it is typically subsampled to obtain approximately uncorrelated samples. 
\\ \\
Although MCMC methods are widely used in applied statistics, their application to large-scale, physics-based Bayesian inverse problems remains very challenging. This is because in each iteration of Algorithm 4 the prior density and likelihood function must be evaluated for the proposed parameters to calculate the acceptance probability \eqref{accept}. The prior density is chosen in advance and is easily evaluated for any parameter values. However evaluating the likelihood function \eqref{17} requires the solution $u(\hat{\bm{x}}_i;\bm{\theta})$ to the PDE or integral equation. This is a key computational bottleneck since a numerical forward solve must be performed at each MCMC iteration. Furthermore a large number of MCMC iterations are required to produce a sufficient an accurate approximation to the posterior. This repeated execution of expensive numerical schemes has so far restricted the MCMC approach to simple models where numerical solvers are not too expensive.
%In our setting the model is \eqref{model}, the data is $\bm{\tilde{z}}=(\bm{\hat{x}},\bm{\hat{z}})$, and the parameters we wish to infer are $\bm{\tilde{\theta}} = \{\bm{\theta},\sigma^2\}$. The likelihood function \eqref{17} in this case requires the solution
\\ \\
To overcome these limitations, we apply the deep surrogate model to Bayesian inference by introducing the surrogate into the Metropolis--Hastings algorithm. Specifically, we approximate $u(\hat{\bm{x}}_i,\bm{\theta})$ by the surrogate $\hat{u}(\hat{\bm{x}}_i,\bm{\theta})$ during the likelihood calculations in the MCMC scheme, thus replacing a numerical solve at each iteration with a cheap evaluation of a known analytic function. Our proposed surrogate model provides a computationally tractable approach to accurately approximate the parametric solutions to PDE and integral equations over the high-dimensional space $\Omega \times \Theta$. Since the surrogates can be evaluated rapidly for any parameter values, they are a powerful tool for Bayesian sampling problems.  Our numerical results in Section 6 demonstrate that this allows for fast posterior sampling when compared to solvers based on grid-based discretisations. Additionally, the separation of the approximation of the parametric solution from the MCMC sampling means that once the surrogate is constructed, it can be stored and used to do inference with different data sets without the need to re-solve the PDE or integral equation. 

\ \\
The next section explores two examples which demonstrate that neural networks allow the accurate and computationally efficient approximation of parametric solutions. The numerical results confirm that this approach significantly accelerates Bayesian inference with an MCMC sampler.

%\subsection{A Note on Point Estimation}
%If point estimates are sufficient then a significant speed-up can be achieved over the Bayesian approach by employing a least squares minimisation over the parameter space to infer the maximum-likelihood parameters assuming Gaussian noise. This technique for making point estimates uses the pre-computed parametric solution, therefore it does not suffer from the issues mentioned in in the previous section. 

\section{Examples}
Here we apply our method to an electrochemistry Bayesian inverse problem that can be represented both as an integral equation and a PDE. We compare our parametric integral equation solver to the deep learning solver for parametric PDEs and an accurate quadrature scheme. Additionally, we apply our methodology to the inference of a function valued parameter used when modelling heat transfer in aircraft compressor turbines. This example demonstrates our approach is capable for high dimensional problems by sampling from the posterior distribution of the parameters of a reduced basis representation of the function of interest. All examples were implemented in TensorFlow using a laptop with a 6 core 3.9GHz CPU and a mobile RTX 2080 GPU.
\subsection{Voltammetry}
Voltammetry is an experimental technique used in electrochemistry to infer the oxidation potential of a chemical. Our description of this problem is based on the construction given in \cite{voltam}. A standard experiment is to apply a potential to an electrode in an electrochemical cell. This causes electron transfer which is measured as an electrical current. This current manifests as the rate of reaction taking place at the electrode surface, which is related to the rate of change in concentration of the chemical species. Using measurements of the current we can deduce information about the chemical system. Specifically, the inverse problem we will consider is to infer the formal potential $E^0$ of the chemical, which is related to the amount of energy that is required to stimulate a reaction. This experiment can be modelled by a PDE or an integral equation, and therefore it allows us to compare the new deep surrogate approach for integral equations to existing PDE-based techniques. 

\subsubsection{PDE Approach}
We first model the system by a 1-dimensional PDE which describes the transport of chemicals through diffusion
\begin{align}
	\frac{\partial C_s}{\partial t}=D_s\frac{\partial ^2C_s}{\partial x^2}.
\end{align}
Here $x \in [0,\infty)$ is the distance from the electrode surface, $C_s(x,t)$ represents the concentration of the chemical species $s$, and $D_s$ is the known scalar diffusion coefficient of species $s$. We seek the solution for all locations $x$ and times $t$ such that $\boldsymbol{x} = (x,t) \in \Omega = [0,\infty)\times[0,\infty)$.
\\ \\
The experiment we will consider involves 2 chemical species, $A$ and $B$. It begins with only chemical $A$ present over the spatial domain at time $t=0$. We assume that the oxidation reaction $A \to B+e^{-}$ (the conversion of $A$ to $B$ through the loss of an electron) takes place at the electrode surface at a rate dependent on the intensity of the applied potential. For simplicity, we assume that $A$ and $B$ share the same constant diffusion coefficient, then after non-dimensionalising appropriately we have a system of equations
\begin{align}
	& & & &\frac{\partial C_a}{\partial t} &= \frac{\partial^2C_a}{\partial x^2}, & (x,t) \in \Omega, & & \label{topde} \\
	& & & &\frac{\partial C_b}{\partial t} &= \frac{\partial^2C_b}{\partial x^2}, 	&		    (x,t) \in \Omega. & &
\end{align}
Here $C_a$ and $C_b$ represent the concentrations of chemicals $A$ and $B$. We assume local conservation of matter, so that
\begin{align}
	&&&& \ \ \ \ \ \ \ \ \ C_a(x,t) + C_b(x,t) = 1, & & (x,t) \in \Omega, &&
\end{align}
therefore it is sufficient to solve for just $C_a(x,t)$. Initially, we assume that only chemical $A$ is present, so the initial condition is
\begin{align}
	&&&& C_a(x,0)&=1, &x \in [0,\infty). && \label{tobc1}
\end{align}
We impose far field boundary conditions for all $t \in [0,\infty)$
\begin{align}
	&&&& C_a(x,t) &\to 1, & x \to \infty. && \label{tobc2}
\end{align}
%
%
%This is imposed through the boundary condition 
%\begin{align}
%    &&&& \ \ \ \ \ \ \ \ \ \frac{\partial C_a}{\partial x} + \frac{\partial C_b}{\partial x} = 0, & & x=0, t\in [0,\infty),    &&
%\end{align}
The final boundary condition at $x=0$ depends on the potential at the electrode. In linear sweep voltammetry a linearly increasing current of the form
\begin{align}
    E(t) = E_{start} + t \label{estart}
\end{align}
is applied. The boundary condition is then
\begin{align}
	&&&& C_a(0,t) &= \frac{1}{1+\text{e}^{E(t)-E^0}},  &t>0. && \label{bc11}
\end{align}
The current is the rate of reaction given by the gradient at the electrode surface
\begin{align}
	I(t)=\frac{\partial C_a}{\partial x}\bigg\rvert_{x=0}. \label{reffff}
\end{align}
The forward problem is to produce this current as a function of the applied potential $E(t)$. The inverse problem is to infer the formal potential $E^0$ of the chemical, given noisy current measurements from experiments.

\subsubsection{Integral Equation Approach}
Linear sweep voltammetry can also be described by an integral equation by taking the Laplace transform of the PDE problem. Letting $\tilde{C_a}(x,s) = \mathcal{L}(C_a(x,t))$ denote the Laplace transform of $C_a(x,t)$ with respect to $t$, it can be shown that
\begin{align}
	\tilde{C_a}(x,s) = c_1(s)\text{e}^{-\sqrt{s}x} + \frac{1}{s},
	\label{49}
\end{align}
for some function $c_1(s)$. To calculate the current we take the Laplace transform of \eqref{reffff}, giving
\begin{align}
	 \tilde{I}(s) = \frac{\partial \tilde{C_a}}{\partial x}\bigg\rvert_{x=0} , \label{5055}
\end{align}
where $\tilde{I}(s)=\mathcal{L}(I(t))$ denotes the Laplace transform of $I(t)$. Differentiating \eqref{49} with respect to $x$, and inserting into $\eqref{5055}$ with $x=0$ we deduce
\begin{align}
	c_1(s) = -\frac{\tilde{I}(s)}{\sqrt{s}},
\end{align}
and therefore
%\begin{align}
%	\tilde{a}(x,s) = -\frac{\tilde{I}(s)}{\sqrt{s}}\text{e}^{-\sqrt{s}x} + \frac{1}{s}.
%\end{align}
%Evaluating this at $x=0$ gives
\begin{align}
	\tilde{C_a}(0,s) = - \frac{\tilde{I}(s)}{\sqrt{s}} + \frac{1}{s}. \label{444444}
\end{align}
To recover $C_a(0,t)$ from \eqref{444444} we apply the inverse Laplace transform to $\tilde{C}_a(0,s)$ using the convolution theorem and identities: $\mathcal{L}^{-1}(\frac{1}{\sqrt{s}})=\frac{1}{\sqrt{\pi t}}$, $\mathcal{L}^{-1}(\frac{1}{s})=1$. This allows us to express the concentration at the electrode as an integral over the current,
%\begin{align}
%	\mathcal{L}^{-1}(\tilde{f}(s)\cdot \tilde{g}(s)) = %\int^t_0f(\tau)g(t-\tau)d\tau
%\end{align}
\begin{align}
	C_a(0,t) = 1 - \frac{1}{\sqrt{\pi}} \int_0^t I(\tau) \cdot \frac{1}{\sqrt{t-\tau}}d\tau.
	\label{n55}
\end{align}
Together with the boundary condition \eqref{bc11} this gives rise to the following Volterra integral equation of the first kind
\begin{align}
	\frac{\sqrt{\pi}}{1+\text{e}^{-(E(t)-E^0)}}=\int_0^t \frac{I(\tau) }{\sqrt{t-\tau}}d\tau. \label{inteqap}
\end{align}
The solution $I(t)$ of this integral equation is the same as the gradient defined in \eqref{reffff} of the solution of the PDE problem (\ref{topde},\ref{tobc1},\ref{tobc2},\ref{bc11}). The objective here is the same as before: we wish to infer the formal potential $E^0$ given measurements of the time-dependent current $I(t)$ from experimentation.

\subsubsection{Results}
The deep surrogate models described in Section 3 and 4 are used to approximate parametric solutions to both the PDE and integral equation formulation of the problems described in Sections 6.1.1 and 6.1.2 respectively. We use $E_{start}=-15$ in \eqref{estart}, and work with the truncated domain $\bar{\Omega} = [0,200] \times [0,25]$. The neural networks were constructed and trained using TensorFlow \cite{tensorflow}. 
\\ \\
For the PDE problem, we applied the methods described in Section 3 to solve the parametric PDE problem using a fully connected neural network with 4 hidden layers of 45 neurons and tanh activation functions. The approximation $\hat{I}(t;E^0)$ to the current defined in \eqref{reffff} is recovered by computing the gradient $\partial_x\hat{C}_a(0,t;E^0)$ at the spatial boundary of the surrogate approximation using automatic differentiation. For the integral equation approach described in Section 6.1.2 we approximate the parametric solution $I(t;E^0)$ of \eqref{inteqap} for $t\in [0,25]$ using the methods described in Section 4. The approximation to the solution and the integrator network were represented by separate fully connected neural networks, each with tanh activation functions and 4 hidden layers containing 45 neurons. The training time for the networks in both cases was under 5 minutes. 

\begin{figure}[H]
\centerline{\includegraphics[width=6.in]{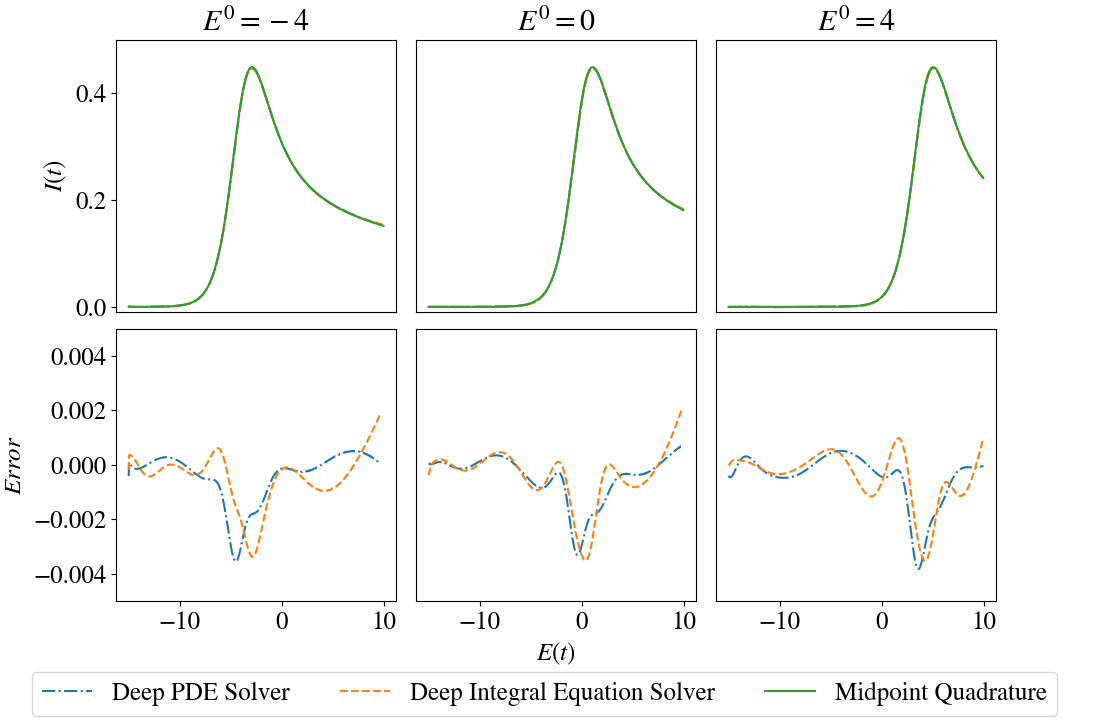}}	
\caption{Comparison of neural network approximations and midpoint quadrature for the voltammetry example in Section 6.1. The top row shows the currents $I(t;E^0)$ approximated by each method for different values of $E^0$. The error in the bottom row is the difference between the neural networks and the reference solution obtained by solving the integral equation with midpoint quadrature.}
\end{figure}
\noindent
As reference solutions $I(t;E^0)$ for each potential $E^0$ we numerically solve the integral equation by approximating the integral in \eqref{inteqap} by a midpoint quadrature rule and solving the resulting system of linear equations. To achieve a high accuracy we use a very small grid spacing in the temporal direction of $\Delta t \approx 5\times 10^{-4}$. These reference solutions are used as reliable approximations of the true solutions and compared to the deep surrogate approach. For $E^0=\{-4,0,4\}$, Figure 1 compares these reference solutions to the parametric solutions achieved by the surrogates. In each case the neural network solutions are close to the references. In particular, we see that the deep learning integral equation solver has comparable accuracy to the established deep learning PDE solver. 
\\ \\
Evaluating both the integral equation and PDE based surrogates is significantly faster than performing quadrature or finite difference at a resolution that achieves the same accuracy. For the integral equation surrogate we are able to perform 178,194 evaluations every second at 100 data locations. In comparison, applying midpoint quadrature with a resolution that achieves the same accuracy as the surrogate we can solve the integral equation problem just 229 times a second using optimised forward substitution methods to solve the lower triangular system. Our PDE surrogate model is slightly slower, achieving 171,943 evaluations at the data points each second, due to requirement to automatically differentiate the approximation to the PDE at the spatial boundary. Solving the 2-dimensional PDE by applying Crank--Nicolson finite differences with sparse linear algebra methods takes 1.51 seconds for a single solve at the same accuracy as the surrogate. Hence in both cases we see a massive a speed-up of 3-6 orders of magnitude using the surrogate models over traditional discretisation schemes.
\\ \\
Having solved the forward problem, we now consider the Bayesian inverse problem. We generate synthetic data for the Bayesian inference task by randomly adding noise to the reference solutions at 100 random time points before sampling. To study the effect of the noise-level of the data on the posterior distributions for $E^0$, Gaussian noise with different variances was used for $E^0 = \{-4,0,4\}$. A uniform $U(-6,6)$ prior distributions was placed on $E^0$, and a $U(0,3)$ prior distribution was placed on $\sigma$. The neural network surrogate solutions were then used in a Metropolis--Hastings sampler \cite{metropolis1953equation} to approximate the joint posterior distribution. Figure 2 shows the simulated data (first row) together with the approximated marginal posterior distributions for $E^0$ obtained using the PDE-based surrogate model (second row). For the integral equation based surrogate the results are presented in the same format in Figure 3. The approximated posterior distributions of $\sigma$ also gave estimates consistent with the true noise levels, however these results are omitted.
\begin{figure}[H]
\centerline{\includegraphics[width=6.in]{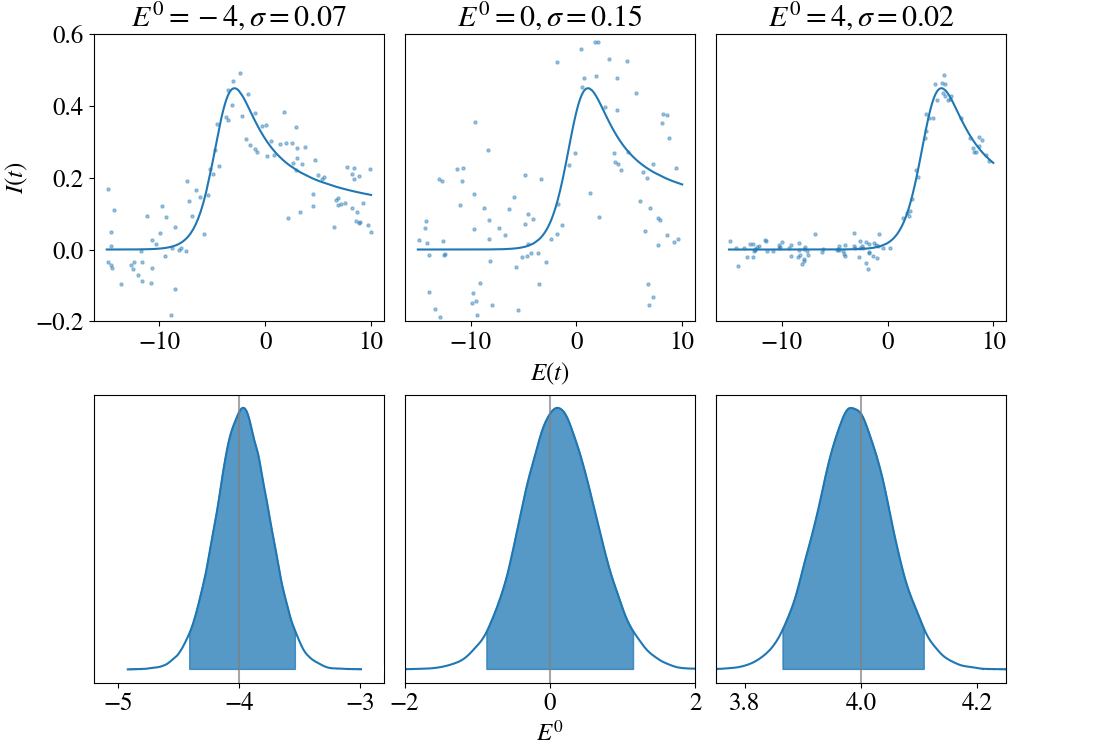}}	
\caption{Top: Simulated data and fitted solutions using the mean of the posterior. The true standard deviations $\sigma$ and oxidation potential $E^0$ used to generate the Gaussian noise are given at the top of each column. Bottom: Approximate posterior distributions and 95\% credible intervals computed using the PDE deep surrogate model. The true values of $E^0$ used to generate the data are marked by vertical lines.}
\end{figure}
\noindent
For both the PDE and integral equation based surrogates, 500,000 Metropolis--Hastings iterations were used to sample from the posterior, and kernel density estimation \cite{kde} was used to visualise the distributions. In each case, the means of the MCMC samples are close to the true parameter used to generate the data, and the width of the credible intervals increases for noisier data as we would expect. The shape and width of the posterior distributions for the integral equation surrogate are consistent with the PDE based surrogate, though we note some slight random variation in their estimates due to the use of different randomly generated noise.
\begin{figure}[H]
\centerline{\includegraphics[width=6.45in]{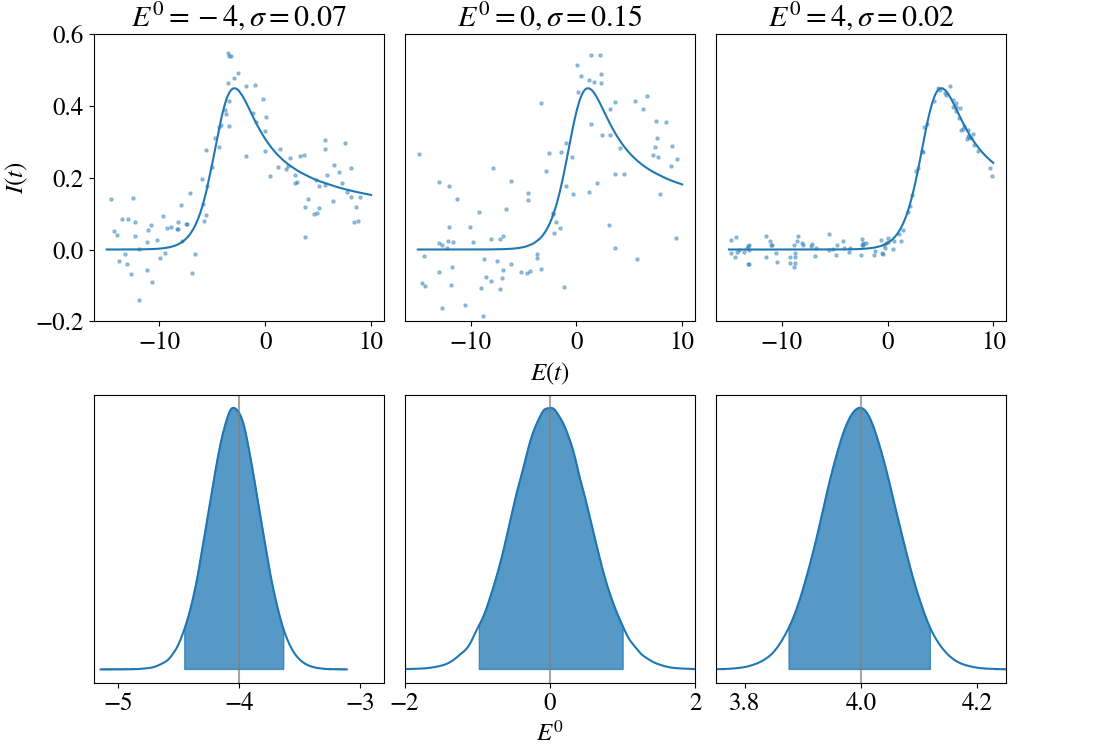}}	
\caption{Top: Simulated data and fitted solutions using the posterior mean. The true oxidation potential $E^0$ and standard deviations $\sigma$ used to generate the Gaussian noise are given at the top of each column. Bottom: Posterior distributions and 95\% credible intervals computed using the integral equation deep surrogate. The true value of $E^0$ used to generate the data are marked by vertical lines. }
\end{figure}

\subsection{Inferring the non-constant Biot number in rotating discs}
Heat transfer in rotating-disc systems encountered in aircraft turbines comprises two mechanisms; heat diffusion in the solid disc and convection of fluid particles. The steady state behaviour of this process can be modelled using the one-dimensional Fin equation \cite{finn}, which only considers variations in the radial direction. The relative effects of diffusion and convection are quantified by a spatially varying parameter known as the \textit{Biot number}. Estimating the Biot number is of interest in aviation applications in order to predict material expansion in aircraft turbines. 
\\ \\
Using suitable units, we consider a rotating disc with inner and outer radii a and 1 respectively, leading to the domain  $\Omega = [a,1]$. In this case, the Fin equation with constant Dirichlet boundary conditions is
\begin{align*}
    \frac{d^2u}{dx^2} + \frac{1}{x}\frac{du}{dx} - Bi\ u &= 0, \qquad x \in \Omega, \\
    u(a) &= u_a, \\
    u(1) &= u_1.\label{biottt} \numberthis
\end{align*}
Here $u=u(x)$ is a non-dimensional temperature, $Bi = Bi(x)$ is the Biot number, $x$ represents radial distance from the centre of rotation, and $u_a,u_1 \in \mathbb{R}$ are fixed Dirichlet boundary conditions. Analytical solutions to \eqref{biottt} are intractable when the Biot number varies over the domain, however it can still be solved numerically in this case. The inverse problem we consider here is to infer the functional parameter $Bi(x)$, given noisy measurements of the temperature profile $u(x)$ at a set of points $x \in \Omega$. 
\\ \\
In principle $Bi(x)$ can be determined numerically if the radial distribution of the temperature u(x) is known, however this is an ill-posed inverse problem where very small uncertainties in the temperature measurements can create large uncertainties in the computed Biot number. This feature of the inverse problem has caused early curve fitting attempts based mean squared error minimisation to infer curves with unphysical behaviour \cite{biot,badcurves}. Similar issues arise if the physics informed deep learning approach for inverse problems \cite{pidl2} is applied directly as Figure 5 below shows. A maximum a posteriori approach (MAP) to the inverse problem of inferring $Bi(x)$ from $u(x)$ indicated that a Bayesian approach can be robust to this sensitivity and allow for reliable inferences \cite{biot}. Although the MAP approach gives us an estimate of $Bi(x)$, it does not reliably quantify uncertainties in this Biot number. To overcome this limitation we apply the deep surrogate method to approximate the posterior distribution of $Bi(x)$.

\subsubsection{Results}
The parameter $Bi(x)$ which is to be inferred lies in the infinite-dimensional space of all continuous function. To make the problem computationally tractable, we consider Biot numbers in the finite dimensional subspace of polynomials of degree up to 15. More specifically, we write the Biot number as a finite sum of monomials 
\begin{align}
    \qquad \qquad \tilde{Bi}(x;\bm{\theta}) = \sum_{n=0}^{15}\theta_nx^n, \qquad x \in \Omega,\ \bm{\theta} \in \Theta.  \label{bipol}
\end{align}
We then use a neural network to generate an approximate parametric solution $\hat{u}(x;\bm{\theta})$ of the PDE, which depends on the position $x\in\Omega=[a,1]$, and the polynomial coefficients $\theta \in \Theta \subset \mathbb{R}^{16}$ that define $\tilde{Bi}(x;\bm{\theta})$. 
\\ \\
Using this representation we approximate the solution $u(x;\bm{\theta})$ of the parametric forward problem by minimising \eqref{loss8}. To represent $\hat{u}(x;\bm{\theta})$ we used a fully connected neural network with a 17-dimensional input layer with arguments $x \in \mathbb{R}$ and $\bm{\theta}\in \mathbb{R}^{16}$, 4 hidden layers with 45 neurons per layer using tanh activation functions. The output layer of the network $\hat{u} (x;\bm{\theta})$ approximates the solution $u(x;\bm{\theta})$. Normal priors were placed on the coefficients of the form $\theta_n \sim N(0,V_n)$, where $V_0 = 20$ and $V_{n} = 20/2^{n-1} \ \ n=1,2,\dots,15$. The decay in the width of the Gaussian priors enforces the smoothness of the functions in the posterior, reducing the ill-posedness of the problem and preventing the nonphysical estimates. The training time for the parametric forward problem over this domain was 14 minutes. Once trained, the surrogate can be evaluated 587304 times each second at 30 data locations. This performance is expected since this network has the same architecture as the surrogate used in the voltammetry integral equation example, however in this example there are 70\% fewer data points leading to an increase in speed of over 3 times. 
\\ \\
Using the trained neural network approximation $\hat{u}(x;\bm{\theta})$ in a Metropolis--Hastings sampler we then sample from the posterior obtained with synthetic data. This data is generated by using finite differences to solve \eqref{biottt} with $Bi(x) = 18e^{x-0.3}$, this non-polynomial is typical of real world data for this scenario (see e.g. \cite{biot}). We add Gaussian noise with standard deviation $\sigma=0.003$ to the solution at 30 equidistant points as shown on the left of Figure 4. Overlaid on this plot is the solution fitted by solving \eqref{biottt} using the $\tilde{Bi}(x;\bm{\theta})$ estimated using the mean of the posterior sample. The right side of Figure 4 shows the estimated $\tilde{Bi}(x;\bm{\theta})$ (blue) and its 95\% credible interval compared to the true $Bi(x)$ (red). The plot demonstrates that the true Biot number lies within the 95\% credible region, and that the fitted solution matches the data well. 
\begin{figure}[H]
\centerline{\includegraphics[width=6.5in]{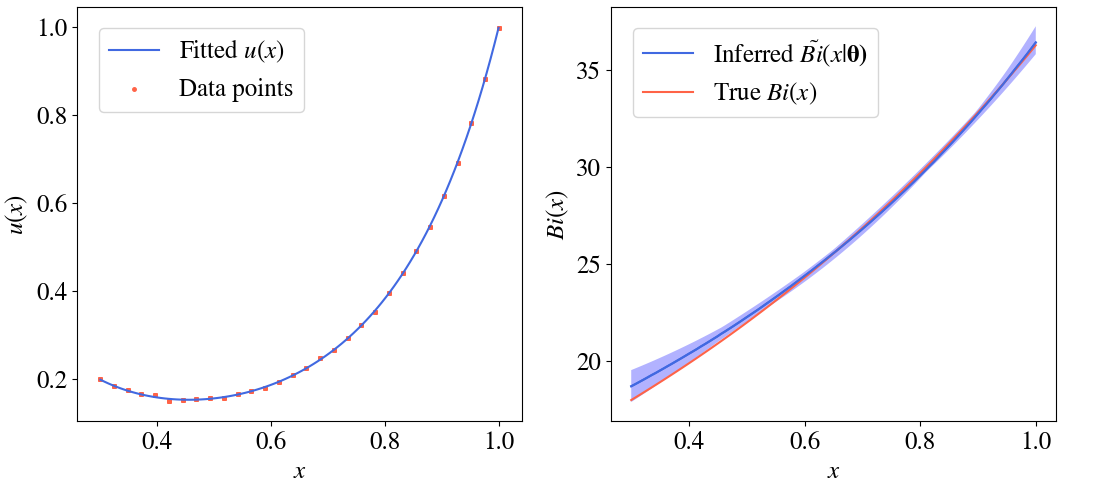}}	
\caption{Left: Fitted heat profile using the inferred Biot number. Right: Inferred $\tilde{Bi}(x;\bm{\theta})$ and 95\% credible interval compared to the true Biot number.}
\end{figure}
\noindent
As a comparison, we implement the physics informed deep learning approach described in \cite{pidl2}. For this we again represent the Biot number as the sum of monomials given in \eqref{bipol} and construct the augmented loss function 
\begin{align*}
            \frac{\nu_1}{N}\sum_{n=1}^N (\mathcal{N}_\theta(\hat{u}_\theta(\bm{x}^n),\bm{x}^n&) - h_\theta(\bm{x}^n))^2 \\
            + \frac{\nu_2}{J}&\sum_{j=1}^J(\hat{u}_\theta(\bm{y}^j) - b_\theta(\bm{y}^j))^2 + \frac{\nu_3}{M}\sum_{i=1}^M(\hat{u}_\theta(\hat{\bm{x}}_i) - \hat{z}_i)^2. \label{lossaug} \numberthis
\end{align*}
The addition of the final term in \eqref{lossaug} encourages the optimiser to minimise the squared difference of the learned solution from the data, this approach has been shown to work well with very large amounts of simulated data. This loss function is minimised with respect to both the parameters of the neural network, and the polynomial coefficients $\bm{\theta}$ simultaneously. Figure 5 shows the estimate achieved by applying this scheme with the same data used in Figure 4. The neural network was trained until convergence, with the weighting coefficients in \eqref{lossaug} set to $\nu_1=\nu_2=\nu_3=1$. 
\begin{figure}[H]
\centerline{\includegraphics[width=6.5in]{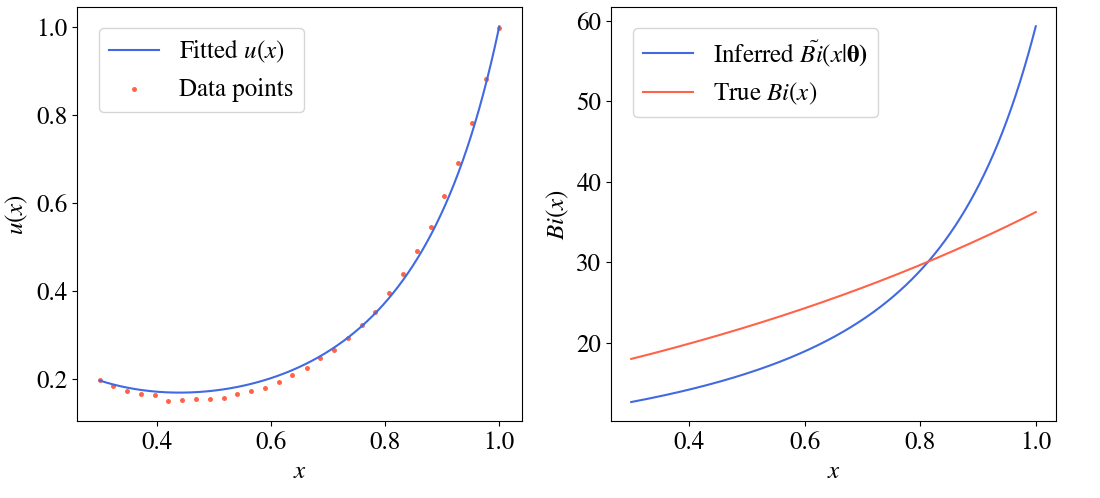}}	
\caption{Right: $\tilde{Bi}(x;\bm{\theta})$ inferred by minimising \eqref{lossaug} compared to true Biot number. Left: Fitted heat profile using the inferred Biot number.} 
\label{poorfit}
\end{figure}
\noindent
Here we see that the fitted solution is close to the data, however the inferred Biot number does not closely match its true value, and would clearly lie outside the 95\% confidence interval computed by our Bayesian approach (see Figure 4). Here we have only 30 noisy data points, and so the minimiser of \eqref{lossaug} must balance the accuracy of the PDE and the interpolation of the data. This risks converging to inaccurate solutions at the expense of better interpolation of the data. In Figure \ref{poorfit} the impact of this behaviour on inference is exacerbated due to the ill-posedness of the inverse problem. The method proposed in \cite{pidl2} also produces no uncertainty estimates over this parameter giving us no indication of confidence in the estimates. This demonstrates that for ill posed inverse problems it is advantageous to compute an accurate surrogate solution first, and that by using Bayesian methodology suitable regularisation and uncertainty quantification can be achieved. 

\section{Conclusion}
In this paper we propose a new deep learning approach for the solutions of integral equations. Our method uses two neural networks, one to approximate the solution of the integral equation, and another to approximate the integral term of these equations. Both networks are trained simultaneously by using mini-batch stochastic gradient descent to minimise a loss function designed such that its minimiser solves the equation. Using this algorithm and existing approaches for PDE problems we described how to construct deep surrogate models to approximate the parametric PDE and integral equation forward problems, by extending the sampling domain to include the parameter space of the problem. Numerical evidence illustrates the accuracy of the parametric integral equation solver when compared to a quadrature based solver. We then applied our surrogate models to sample from the posterior distribution of PDE and integral equation based Bayesian inverse problems. We demonstrated that our surrogate models can be evaluated rapidly, achieving a sampling acceleration of several orders of magnitude at the same accuracy when compared to more traditional discretisation schemes. The examples that we considered in this paper show that empirical posterior distributions of the model parameters achieved by applying this scheme are consistent with the true parameters.
\\ \\
In contrast to traditional grid-based methods, for which the computational complexity grows exponentially with the dimension, the deep surrogate model is efficient at approximating solutions to high dimensional parametric forward problems. Additionally, since our method is trained to satisfy the equations directly over the entire parameter space, it approximates the solution more accurately than surrogate models based on interpolation. Our framework is flexible enough to be applied to a wide range of PDEs and integral equations, and extensions to higher dimension integral equations, problems with Neumann boundary conditions, or systems of equations, are possible using the principles described. Once training is complete the deep surrogate model can be evaluated rapidly, and in parallel, making deep surrogates ideal candidates to accelerate MCMC schemes for Bayesian inverse problems. 
\\ \\
Upcoming research on this topic will investigate how deep surrogate models can be applied in more complex and varied scenarios. A key benefit of the surrogate over traditional methods is its differentiability with respect to the PDE or integral equation model parameters. Future work will utilise this to develop bespoke samplers using gradient based transition kernels and delayed acceptance criteria that mathematically guarantees the accuracy of the posterior samples in large-scale problems. Other possible extension include adaptive surrogate construction and the interaction of physics (represented by deep surrogate models) with applied statistical problems in environmental sciences. 

\subsection*{Acknowledgements}
Teo Deveney  is supported by a scholarship from the EPSRC Centre for Doctoral Training in Statistical
Applied Mathematics at Bath (SAMBa), under the project EP/L015684/1.

\bibliography{main2}
\bibliographystyle{plain}
\end{document}